\documentclass[graybox]{svmult}

\usepackage{sidecap}
\usepackage{makeidx}         
\usepackage{graphicx}        
\usepackage{multicol}        
\usepackage[bottom]{footmisc}

\usepackage{newtxtext}       %
\usepackage{amsmath} 
\usepackage{amsfonts} 

\makeindex             

\begin{document}

\title*{Discrete analysis of Schwarz Waveform Relaxation for a simplified air-sea coupling problem with nonlinear transmission conditions}
\titlerunning{SWR for a simplified air-sea coupling problem}
\author{S. Clement \and F. Lemarié \and E. Blayo}
\institute{S. Clement \at S. Clement, Univ Grenoble Alpes, CNRS, Inria, Grenoble INP, LJK, Grenoble, France, \email{simon.clement@grenoble-inp.org}
\and F. Lemarié\at F. Lemarié, Univ Grenoble Alpes, Inria, CNRS, Grenoble INP, LJK, Grenoble, France, \email{florian.lemarie@inria.fr}
\and E. Blayo \at E. Blayo, Univ Grenoble Alpes, CNRS, Inria, Grenoble INP, LJK, Grenoble, France, \email{eric.blayo@univ-grenoble-alpes.fr}}
%
%
\maketitle

\abstract*{
In this study we present a non-overlapping Schwarz waveform relaxation (SWR) method
applied to a one dimensional model problem representative of the coupling between 
the ocean and the atmosphere. This problem includes nonlinear interface conditions 
analogous to a quadratic friction law. We study the convergence of the corresponding 
SWR at a semi-discrete level for a linear friction and for a linearized quadratic
friction at the interface. Using numerical experiments we show that the convergence 
properties in the linearized quadratic friction case are very close to the ones 
obtained with the full nonlinear problem for the range of parameter values of 
interest. We investigate the possibility to improve the convergence speed by 
adding a relaxation parameter at the interface.}

\section{Introduction}
Schwarz-like domain decomposition methods are very popular in 
mathematics, computational sciences and engineering notably for 
the implementation of coupling strategies. Such an iterative method
has been recently applied in a state-of-the-art Earth 
System Model (ESM) to evaluate the consequences of inaccuracies 
in the usual ad-hoc ocean-atmosphere coupling algorithms 
used in realistic models \cite{clement_mini_02_Marti_etal_2020}. For such 
a complex application it is challenging to have an a priori 
knowledge of the convergence properties of the Schwarz method. 
Indeed coupled problems arising in ESMs often 
exhibit sharp turbulent boundary layers whose parameterizations 
lead to peculiar transmission conditions. The objective in
this paper is to study a model problem representative 
of the coupling between the ocean and the atmosphere, including discretization 
and so-called bulk interface conditions which are analogous to a 
quadratic friction law. Such a model is introduced in Sec. \ref{clement_mini_02_sec:couplingProblem}
and its discretization, as done in state-of-the-art ESMs, 
is described in Sec. \ref{clement_mini_02_sec:discreteProblem}. In the semi-discrete
case in space we conduct in Sec. \ref{clement_mini_02_sec:conv-lin} a convergence analysis 
of the model problem first with a linear friction and then with a 
quadratic friction linearized around equilibrium solutions. 
Finally, in Sec. \ref{clement_mini_02_sec:num-exp}, numerical experiments in the 
linear and nonlinear case are performed to illustrate the relevance 
of our analysis. \par
\section{Model problem for ocean-atmosphere coupling}\label{clement_mini_02_sec:couplingProblem}
We focus on the dynamical part of the oceanic and atmospheric primitive equations and neglect the horizontal variations of 
the velocity field, which leads to a model problem  depending 
on the vertical direction only. This assumption, commonly made 
to study turbulent mixing in the boundary layers near the 
air-sea interface, is justified because of the large disparity 
between the vertical and the horizontal spatial scales in these layers. 
We consider the following diffusion problem accounting for Earth’s 
rotation ($f$ is the Coriolis frequency and $\mathbf{k}$ a vertical unit vector):
\begin{equation*}
\left\{
\begin{array}{rcll}
\partial_t \mathbf{u} + f \mathbf{k} \times \mathbf{u} - \partial_z \left( \nu(z,t) \partial_z \mathbf{u} \right) &=& \mathbf{g}, & \mbox{in}\;\Omega \times (0,T), \\
\mathbf{u}(z,0) &=& \mathbf{u}_0(z), & \forall z \; \mbox{in}\; \Omega, \\
\mathbf{u}(H_o,t) &=& \mathbf{u}_o^\infty(t), ~ \mathbf{u}(H_a,t) = \mathbf{u}_a^\infty(t), ~ & t \in (0,T),  
\end{array}
\right.
\end{equation*}
with $\mathbf{u} = (u,v)$ the horizontal velocity vector,
$\nu(z,t) > 0$ the turbulent viscosity and 
$\Omega = (H_o,H_a)$ a bounded open subset of $\mathbb{R}$
containing the air-sea interface $\Gamma = \{z = 0\}$.
In the ocean and the atmosphere, which are  turbulent fluids,
the velocity field varies considerably in the few 
meters close to the interface (in a region called 
\textit{surface layer}). The cost of an
explicit representation of the surface layer in numerical simulations being 
unaffordable, this region is numerically accounted for 
using wall laws a.k.a. log laws (e.g. \cite{clement_mini_02_Mohammadi_etal_1998}). This approach, 
traditionally used to deal with solid walls, is also used
in the ocean-atmosphere context, with additional complexity 
arising from the stratification effects \cite{clement_mini_02_Pelletier_etal_2021}.
In this context wall laws are referred to as \textit{surface layer} 
parameterizations. The role of such parameterizations is to provide 
$\nu \partial_z \mathbf{u}$ on the upper and lower interfaces of the 
surface layer as a function of the difference of fluid velocities.
Thus the coupling problem of interest should be understood as a domain
decomposition with three non-overlapping subdomains. For the sake 
of convenience the velocity vector $\mathbf{u} = (u,v)$ 
is rewritten as a complex variable $U = u + i v$. Then the model problem reads   
\begin{equation}
\label{clement_mini_02_eq:cplProblem}
\begin{array}{rcll}
\partial_t U_j + i f U_j - \partial_z \left( \nu_j(z,t) \partial_z U_j \right) &=& g_j,
~~~~~~~~~~ (j=o,a)
& \mbox{in}\;\Omega_j \times (0,T)\\
U_j(H_j,t) &=& U_j^\infty(t),  & t \in (0,T), \\ 
U_j(z,0) &=& U_0(z), & \forall z \; \mbox{in}\; \Omega_j, \\
\rho_o \nu_o \partial_z U_o(\delta_o,t) = \rho_a \nu_a \partial_z U_a(\delta_a,t)
&=& {\cal F}_{\rm sl}( U_a(\delta_a,t)-U_o(\delta_o,t) ), \;\; & t \in (0,T)
\end{array}
\end{equation}
where $\Omega_o= (H_o,\delta_o)$, $\Omega_a = (\delta_a,H_a)$, and 
${\cal F}_{\rm sl}$ is a parameterization function for the surface layer 
extending over $\Omega_{\rm sl} = (\delta_o,\delta_a)$. A typical 
formulation for ${\cal F}_{\rm sl}$ is
\[
{\cal F}_{\rm sl}( U_a(\delta_a,t)-U_o(\delta_o,t) ) = 
\rho_a C_D | U_a(\delta_a,t)-U_o(\delta_o,t) | (U_a(\delta_a,t)-U_o(\delta_o,t))
\]
which corresponds to a quadratic friction law with $C_D$ 
a drag coefficient (assumed constant in the present study).
%
%
Geostrophic winds and currents are used in this study as
source terms and boundary conditions. 
Geostrophic equilibrium is the stationary state for which the 
Coriolis force compensates for the effects of gravity.
It corresponds to the large scale dynamics of ocean and 
atmosphere, and leads to reasonable values of the solution $U$. \par
The well-posedness of \eqref{clement_mini_02_eq:cplProblem} has been studied in 
\cite{clement_mini_02_Thery_2021}
where it is proved that its stationary version admits a unique solution for realistic values of the parameters. The study of the nonstationary case is much more challenging: numerical experiments tend to confirm this well-posedness, but  with no theoretical proof.

\section{Discretized coupled problem}\label{clement_mini_02_sec:discreteProblem}
\subsection{Implementation of the surface layer}
\begin{SCfigure}
    \centering
    \includegraphics[scale=0.76]{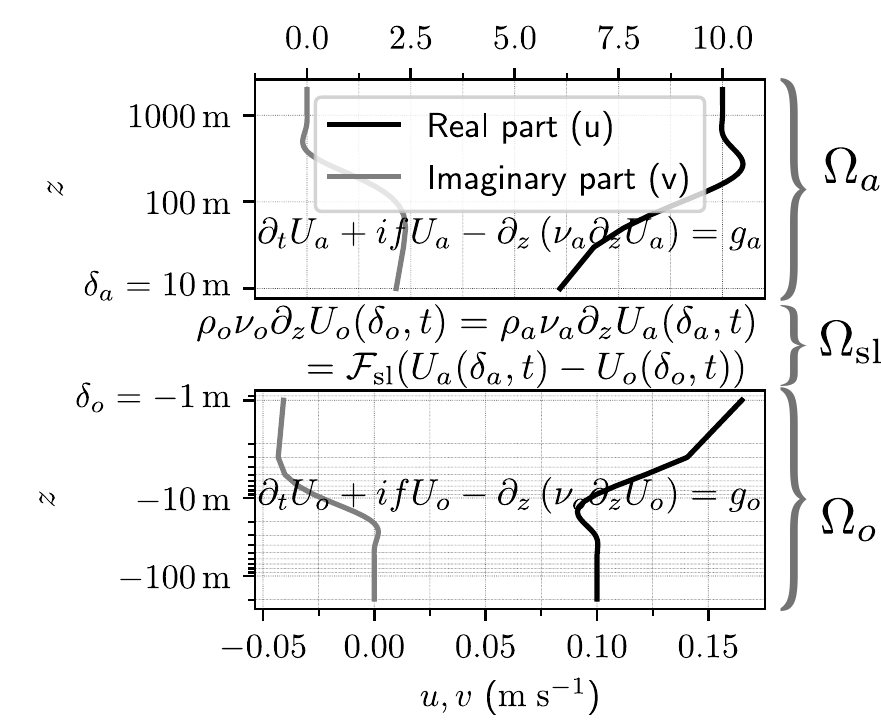}
    \caption{Discrete representation of the three domains $\Omega_a, \Omega_{\rm sl}, \Omega_o$ together with a typical stationary state. Note the different scales for $(u,v)$ in the ocean and in the atmosphere.}
    \label{clement_mini_02_fig:presentationDomains}
\end{SCfigure}

As described in Sec. \ref{clement_mini_02_sec:couplingProblem}, the full domain 
$\Omega$ is split into three parts: $\Omega_o$ in the ocean, 
$\Omega_a$ in the atmosphere and $\Omega_{\rm sl}$ a thin domain 
containing the interface (see Fig. \ref{clement_mini_02_fig:presentationDomains}).
The role of  
$\Omega_{\rm sl}$ is to provide $\rho_j \nu_j \partial_z U_j$ 
at $z=\delta_j$ ($j=o,a$) as a function of fluid velocities at the 
same locations. However, in state-of-the-art climate models, the discretization 
is based on an approximate form of the coupled 
problem \eqref{clement_mini_02_eq:cplProblem}. For practical reasons, the computational
domains are $\widetilde{\Omega}_o = (H_o,0)=\Omega_o \bigcup (\delta_o, 0)$ and 
$\widetilde{\Omega}_a = (0,H_a) = (0,\delta_a ) \bigcup \Omega_a$, and the locations 
of the lower and upper boundaries of the surface layer ($z=\delta_j$)   
are assimilated to the centers of the first grid cells (i.e. $\delta_o=-h_{\rm o}/2$ 
and $\delta_a=h_{\rm a}/2$ with $h_{\rm o}$ and $h_{\rm a}$ the thicknesses of the first grid cell in each 
subdomain), where the  values of the velocity closest to the interface are available. 
%
%
%
%
Typical resolutions in the  models are 
$\delta_a=h_{\rm a}/2=10\;{\rm m}$ and $\delta_o=-h_{\rm o}/2=-1\;{\rm m}$.
At a discrete level, the transmission condition in \eqref{clement_mini_02_eq:cplProblem}
is replaced by
\begin{equation}\label{clement_mini_02_eq:discreteCouplingConditions}
    \rho_{\rm o} \nu_{\rm o} \partial_z U_{\rm o}(0,t)
    = \rho_{\rm a} \nu_{\rm a} \partial_z U_{\rm a}(0,t)
    = \rho_{\rm a} \alpha \left(
    U_{\rm a}\left( \frac{h_{\rm a}}{2},t \right) - U_{\rm o}\left(-\frac{h_{\rm o}}{2},t \right)
    \right)
\end{equation}
where 
$\alpha= C_D \left| U_{\rm a}\left( \frac{h_{\rm a}}{2},t \right) - U_{\rm o}\left(-\frac{h_{\rm o}}{2},t \right) \right|$ for the \textit{nonlinear} case.
In the following, for the analysis in Sec. \ref{clement_mini_02_sec:conv-lin}, we  consider 
a \textit{linear} friction where $\alpha$ is assumed constant and a 
quadratic friction \textit{linearized} around equilibrium solutions.
\subsection{Schwarz Waveform Relaxation}
As discussed for example in \cite{clement_mini_02_Marti_etal_2020},
current ocean-atmosphere coupling methods can actually be seen as a 
single iteration of a Schwarz Waveform Relaxation (SWR) algorithm. 
SWR applied to the coupling problem presented 
in Sec. \ref{clement_mini_02_sec:couplingProblem} with the transmission conditions 
 \eqref{clement_mini_02_eq:discreteCouplingConditions} and constant viscosity in each subdomain 
 reads:
%
\begin{subequations}
\label{clement_mini_02_eq:SWRbulk}
\begin{align}
(\partial_t + if) &U^k_j - \nu_j\partial_z \phi^k_j  = g_j, 
\hspace{3.4cm} \mbox{in}\;\widetilde{\Omega}_j \times (0,T) \label{clement_mini_02_eq:ref} \\
U^k_j(z,0) &= U_0(z),   \hspace{4.8cm}  \forall z \in \widetilde{\Omega}_j  \\
U^k_j(H_j, t) &= U^\infty_j, \hspace{5.1cm}  t \in [0,T]\\
\nu_{\rm a}\phi_{\rm a}^{k}(0,t) &=  \alpha^{k-1} 
\left( U_{\rm a}^{k-1+\theta}\left(\frac{h_{\rm a}}{2},t\right) - U_{\rm o}^{k-1}\left(-\frac{h_{\rm o}}{2},t\right)\right), 
\hspace{0.12cm} t \in [0,T] \label{clement_mini_02_eq:bulkInterfaceCondition} \\
\rho_{\rm o} \nu_{\rm o}\phi^k_{\rm o}(0,t) &= \rho_{\rm a}
\nu_{\rm a}\phi^k_{\rm a}(0,t), \hspace{3.98cm} t \in [0,T] \label{clement_mini_02_eq:fluxInterfaceCondition}
\end{align}
\end{subequations}
where $j={\rm a},{\rm o}$, $\phi_j = \partial_z U_j$, and 
$U_{\rm a}^{k-1+\theta} = \theta U_{\rm a}^{k} + (1-\theta) U_{\rm a}^{k-1}$
with $\theta$ a relaxation parameter (interpolation for $0\leq \theta \leq 1$ or extrapolation for $\theta>1$). 
%
%
At each iteration, \eqref{clement_mini_02_eq:fluxInterfaceCondition} ensures that the 
kinetic energy is conserved at the machine precision in the coupled 
system which is a major constraint for climate models. 
In \eqref{clement_mini_02_eq:bulkInterfaceCondition},
the presence of the parameter $\theta$ makes it resemble to a Dirichlet-Neumann Waveform Relaxation algorithm.
Indeed, if \eqref{clement_mini_02_eq:bulkInterfaceCondition} is replaced
by $U_{\rm a}^k = \theta U_{\rm o}^{k-1} + (1-\theta) U_{\rm a}^{k-1}$ the DNWR algorithm is retrieved, as examined in the continuous case in \cite{clement_mini_02_kwok2013} and in the discrete case in
\cite{clement_mini_02_mongeMultirate2021}.
However \eqref{clement_mini_02_eq:bulkInterfaceCondition}
involves both $\phi_{\rm a}^k$ and 
$U_{\rm a}^{k-1+\theta}$:
the $\theta$ parameter appears thus here within (close to Robin) condition  
($\nu_{\rm a} \phi_{\rm a}(0) - \alpha \theta U_{\rm a}(h_{\rm a}/2) = \ldots$),
i.e. the relaxation is not performed directly on the converging variable
which leads to convergence properties different from the DNWR case, 
as shown in Sec. \ref{clement_mini_02_sec:conv-lin}.\par
In the following, centered finite difference schemes in space are used with constant space steps $h_j$.
Derivatives are $\phi_j(z,t) = \frac{U_j(z+h_j/2,t) - U_j(z-h_j/2, t)}{h_j}$
and the semi-discrete version of \eqref{clement_mini_02_eq:ref} in the homogeneous case is
\begin{equation}\label{clement_mini_02_eq:spaceTimeScheme}
(\partial_t + if) U_j (z, t) = \nu_j \frac{\phi_j(z+h_j/2,t) - \phi_j(z-h_j/2,t)}{h_j}
\end{equation}
\section{Convergence analysis} \label{clement_mini_02_sec:conv-lin}
In this section we conduct a convergence analysis of the SWR algorithm 
\eqref{clement_mini_02_eq:SWRbulk} first with $\alpha$ a constant and then in a more 
complicated case where the problem is linearized around its 
equilibrium solutions. In the following we systematically make the assumption 
that the space domain is of infinite size (i.e. $H_j\to\infty$) for the 
sake of simplicity.
\\[3mm]
\textbf{Linear friction case ($\alpha = {\rm const}$)}\hspace*{5mm} 
%
%
We assume in this paragraph that $\alpha=\alpha_{c}$ with $\alpha_{c}$ 
a constant independent of $U_j$ and we study the system satisfied by the errors
(i.e. $g_j, U_0, U^\infty =0$). 
The Fourier transform in time of the finite difference scheme \eqref{clement_mini_02_eq:spaceTimeScheme} 
yields  
$\widehat{U}_{\rm a}(h_{\rm a}/2) = \nu_{\rm a} \frac{\widehat{\phi}_{\rm a}(h_{\rm a}) - \widehat{\phi}_{\rm a}(0)}{i(\omega+f) h_{\rm a}}$
with $\omega \in \mathbb{R}$ the frequency variable. 
After simple algebra, the transmission condition 
\eqref{clement_mini_02_eq:bulkInterfaceCondition} in Fourier space expressed 
in terms of the $\widehat{\phi}_j$ is 
\begin{equation} \label{clement_mini_02_eq:discreteBulkFourier}
\begin{aligned}
 \left( \frac{\chi_a\nu_a}{h_a}+\theta \alpha_c \right) \widehat{\phi}^k_{\rm a}(0) - \theta \alpha_c
\widehat{\phi}^k_{\rm a}(h_a) = &(1-\theta) \alpha_c 
(\widehat{\phi}^{k-1}_{\rm a}(h_{\rm a}) - \widehat{\phi}^{k-1}_{\rm a}(0)) \\
&- \alpha_c \frac{h_a \nu_o}{h_o \nu_a} (\widehat{\phi}^{k-1}_{\rm o}(0) - \widehat{\phi}^{k-1}_{\rm o}(-h_{\rm o}))
\end{aligned}
\end{equation}
with $\chi_j=\frac{i (\omega+f) h_j^2}{\nu_j}$.
%
%
A discrete analysis of the finite difference scheme \eqref{clement_mini_02_eq:spaceTimeScheme} in the
frequency domain (e.g. \cite{clement_mini_02_Wu2017}) leads to 
$\widehat{\phi}_{\rm o}^k(-m h_{\rm o}) = A_k (\lambda_{\rm o}+1)^m$
and
$\widehat{\phi}_{\rm a}^k(m h_{\rm a}) = B_k (\lambda_{\rm a}+1)^m$ 
with $\lambda_j = \frac{1}{2}\left(\chi_j - \sqrt{\chi_j} \sqrt{\chi_j + 4}\right)$ 
and $m$ the space index. The convergence factor of SWR is then
the rate at which $A_k$ or $B_k$ tends to 0.
Combining \eqref{clement_mini_02_eq:discreteBulkFourier} with
the Fourier transform in time of \eqref{clement_mini_02_eq:fluxInterfaceCondition},
we get the evolution of $B_k$ which eventually leads to 
the following convergence factor: 
%
\begin{equation} \label{clement_mini_02_eq:cvBk}
\xi = \left|\frac{B_k}{B_{k-1}}\right| = 
\left|\frac{
\left(1-\theta\right)
+ \epsilon\frac{h_{\rm a} \lambda_{\rm o}}{h_{\rm o} \lambda_{\rm a}}
}{
\frac{\nu_{\rm a}\chi_{\rm a}}{\alpha_c h_{\rm a} \lambda_{\rm a}} - \theta}\right|,
\end{equation}
where $\epsilon = \frac{\rho_{\rm a}}{\rho_{\rm o}} \approx 10^{-3}$
in the ocean-atmosphere  context.
Note that the convergence factor \eqref{clement_mini_02_eq:cvBk} differs significantly from the semi-discrete convergence factor 
$\xi_{\rm DNWR} = \left|1 - \theta_{\rm DNWR} \left(
1 - \epsilon h_{\rm a} \lambda_{\rm o} / (\lambda_{\rm a}h_{\rm o})\right) \right|$ of 
the DNWR algorithm. 
Moreover, it can be found that 
\[
 \underset{(\omega+f) \rightarrow 0}{\mathrm{lim}} \xi = \frac{1}{\theta} \left| 1-\theta+\epsilon \sqrt{\frac{\nu_{\rm a}}{\nu_{\rm o}}}  \right| = \xi_0, \qquad
 \underset{(\omega+f) \rightarrow \infty}{\mathrm{lim}} \xi = 0.
\]
As $\omega + f \to 0$ the asymptotic value 
$\xi_0$ depends on $\theta$: it is $+\infty$ for $\theta=0$ (i.e. a fast divergence), and $\xi_0 = \epsilon\sqrt{\frac{\nu_{\rm a}}{\nu_{\rm o}}}$ for $\theta=1$.
When $\omega\to \infty$, the convergence factor tends to zero 
(i.e. the convergence is fast for high frequencies).
Whatever $\omega$, it can be shown that the value $\xi_0$ is an upper bound 
of the convergence factor when $\theta \leq 1$
if $\sqrt{\frac{\nu_{\rm o}}{\nu_{\rm a}}} \leq \frac{h_{\rm o}}{h_{\rm a}}$,
the latter condition being easily satisfied.
Since we have $\epsilon \approx 10^{-3}$, 
the convergence is fast for $\theta = 1$ whereas $\epsilon$
does not play any role for $\theta = 0$.
%
The optimal parameter $\theta_{\rm opt}$ for low frequencies is $1 +  \epsilon\sqrt{\frac{\nu_{\rm a}}{\nu_{\rm o}}}$ which is very close to $1$.
\\[2mm]
\textbf{Linearized quadratic friction case}\hspace*{5mm}
The analysis of the nonlinear quadratic friction case (i.e. with $\alpha = C_D\left|\right.U_{\rm a}\left(h_{\rm a} / 2, t\right) - U_{\rm o}\left(- h_{\rm o}/2, t\right)\left.\right|$)
cannot be pursued through a Fourier transform.
We thus consider the linearization of the problem 
around a stationary state $U^e_j, \phi^e_j$ satisfying \eqref{clement_mini_02_eq:cplProblem}:
assuming that $U^k_j (\pm h_j/2, t)$ is in a neighborhood of $U^e(\pm h_j/2)$, the modulus in $\alpha$
is non-zero and
we can differentiate $\alpha$.
Differences with the stationary state are noted $\delta \phi_j^k = \phi^k_j(0,t) - \phi_j^e(0)$ and $\delta U_j^k = U_j^k(\pm h_j/2, t) - U^e_j(\pm h_j/2)$.
After some algebra, the linearized transmission operator reads
\begin{equation}
\label{clement_mini_02_eq:linearised}
\begin{aligned}
\nu_{\rm a}\delta\phi_{\rm a}^k =  \alpha^e 
\left( \left(\frac{3}{2} - \theta \right)\delta U_{\rm a}^{k-1} \right.
&+ \theta\, \delta U_{\rm a}^{k}
- \frac{3}{2} \delta U_{\rm o}^{k-1}
\\
&\left.+ \frac{1}{2}\frac{U_{\rm a}^{{\color{black} e}} - U_{\rm o}^{{\color{black} e}}}{\overline{U_{\rm a}^{{\color{black} e}} - U_{\rm o}^{{\color{black} e}}}}\,
\overline{
\delta U_{\rm a}^{k-1} - \delta U_{\rm o}^{k-1}}
\right)
\end{aligned}
\end{equation}
with $\alpha^e=C_D\left|U_{\rm a}^e(h_{\rm a}/2) - U_{\rm o}^e(-h_{\rm o}/2)\right|$.
Following the derivation in the previous paragraph,
we find that the convergence factor $\xi^q$ in the linearized 
quadratic friction case differs from one iteration to another 
(it is indeed a function of $\frac{\overline{B_{k-1}(-\omega)}}{B_{k-1}(\omega)}$).
However, for $(\omega+f) \rightarrow 0$ the term 
$ \frac{1}{2}\frac{U_{\rm a}^e - U_{\rm o}^e}{\overline{U_{\rm a}^e - U_{\rm o}^e}}\,
\overline{
\delta U_{\rm a}^{k-1} - \delta U_{\rm o}^{k-1}}$
vanishes, therefore the asymptotic convergence 
rate $\xi_0^q$ is independent of the iterate:
\[
 \underset{(\omega+f) \rightarrow 0}{\mathrm{lim}} \xi^q = \frac{1}{\theta} \left| \frac{3}{2}-\theta+\frac{3}{2}\epsilon\sqrt{\frac{\nu_{\rm a}}{\nu_{\rm o}}} \right| = \xi_0^q, \qquad
 \underset{(\omega+f) \rightarrow \infty}{\mathrm{lim}} \xi^q = 0.
\]
The convergence is fast for high frequencies,
as in the linear friction case. However
the optimal parameter for $(\omega+f) \rightarrow 0$ is
here $\theta_{\rm opt}^q = \frac{3}{2} + \frac{3}{2}\epsilon\sqrt{\frac{\nu_{\rm a}}{\nu_{\rm o}}} $.
It is different from the optimal parameter $\theta_{\rm opt}$ obtained with linear friction: for typical values of the ocean-atmosphere coupling problem, $\theta_{\rm opt}^q$ is close to $\frac{3}{2}$. The asymptotic value $\xi_0^q$ is not an upper bound of the convergence factor but it is  a good choice 
for $\theta_{\rm opt}^q$.

\section{Numerical experiments}
\label{clement_mini_02_sec:num-exp}
\begin{figure}
    \centering
    \includegraphics[scale=0.65]{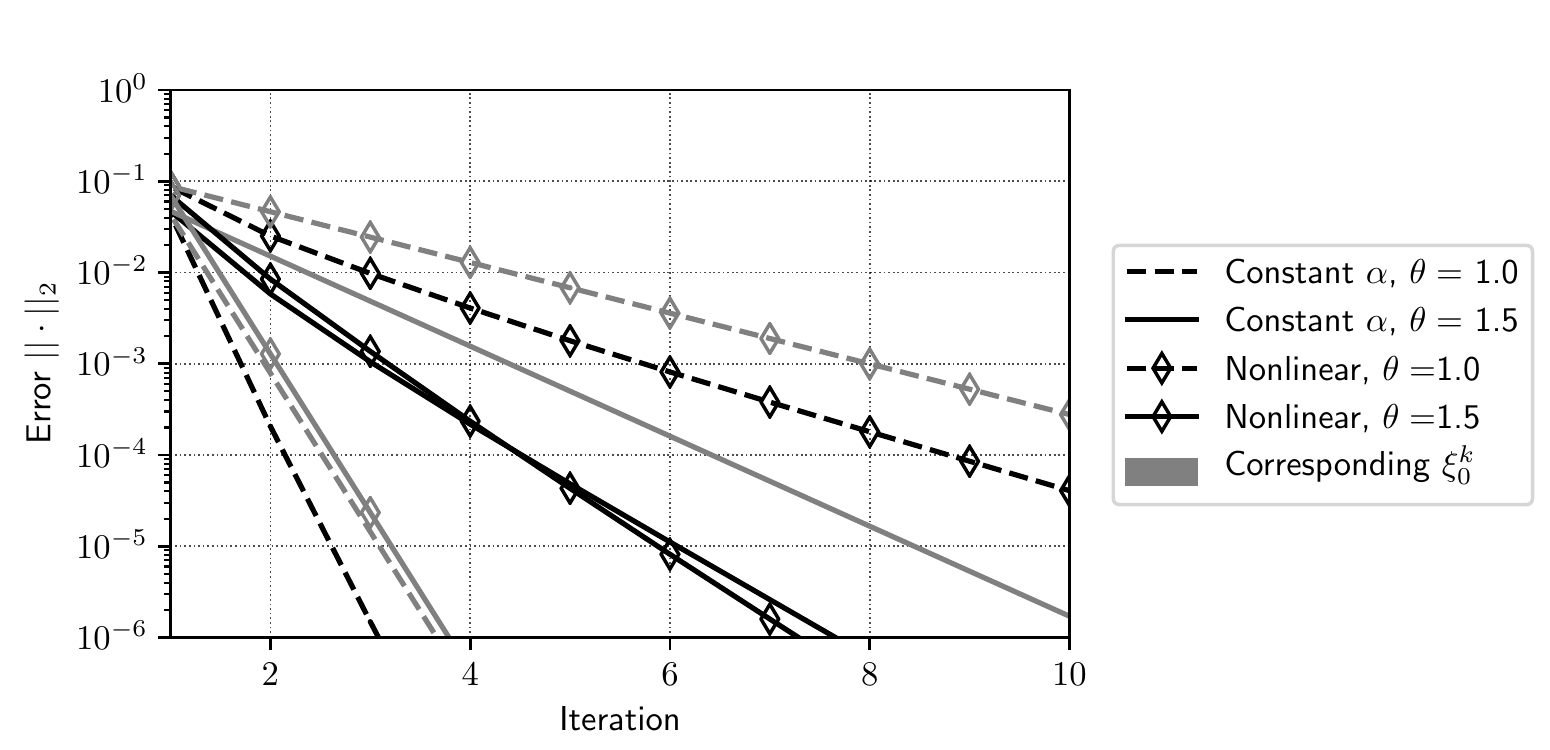}
    \caption{Evolution of the  $L^2$ norm of the errors. Black lines represent the observed convergence; grey lines are the estimated convergence with slopes $\xi_0$ for linear cases and $\xi_0^q$ for quadratic cases.}
    \label{clement_mini_02_fig:evolutionErrNonlinear}
\end{figure}

\begin{figure}
    \centering
    \includegraphics[scale=0.65]{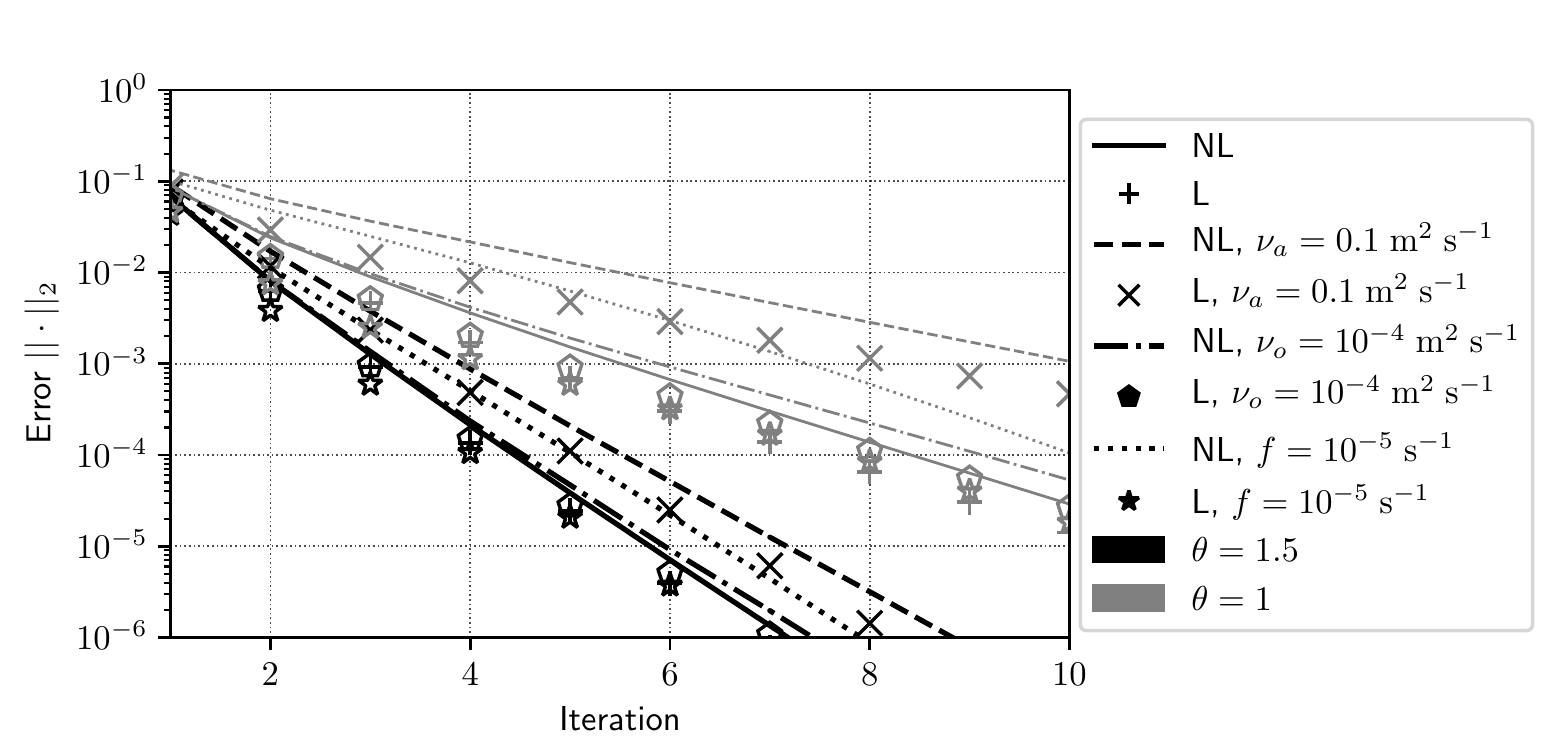}
    \caption{Evolution of the  $L^2$ norm of the errors with linearized (L) and nonlinear (NL) transmission conditions. The legend indicates the changes in the parameters for each case. }
    \label{clement_mini_02_fig:robustesseEvolutionErrNonlinear}
\end{figure}

The aim of this section is to illustrate the influence of the parameter $\theta$, in the linear and quadratic friction cases. 
The stationary state $U_j^e$ is used to compute $\alpha_c = \alpha^e = C_D|U_{\rm a}^e(\frac{h_{\rm a}}{2}) - U_{\rm o}^e (\frac{h_{\rm o}}{2})|$ in the linear case.
Parameters of the problem are taken as realistic: $C_D = 1.2\times 10^{-3}$, the space steps are $\frac{h_{\rm a}}{2} = 10\; {\rm m}$, $\frac{h_{\rm o}}{2} = 1 \;{\rm m}$,
the time step  is $60 \;{\rm s}$,
the size of the time window $T$ is 1 day ($1440\Delta t$) and the computational domains sizes are $H_o = H_a = 2000 \; {\rm m}$
(100 and 1000 nodes respectively in $\Omega_a$ and $\Omega_o$).
The Coriolis parameter is 
$f=10^{-4}\; {\rm s^{-1}}$ and the diffusivities are 
$\nu_{\rm a} = 1\; {\rm m^2}\;{\rm s}^{-1}, 
\nu_{\rm o} = 3\times 10^{-3} \;{\rm m^2}\;{\rm s}^{-1}$. 
 $U_j^\infty$
are  set to constant values of $10\;{\rm m}\;{\rm s}^{-1}$ in the 
atmosphere and $0.1\;{\rm m}\;{\rm s}^{-1}$ in the ocean, while the forcing terms 
$g_j=i f U^\infty_j$ and the initial condition $U_0(z)=U_j^e(z)$.
SWR is 
initialized at the interface with a white noise around the interface value
of the initial condition.
Figure \ref{clement_mini_02_fig:evolutionErrNonlinear} shows the evolution of the error for two choices of $\theta$. The theoretical
convergence according to $\xi_0$ is
also displayed: $\sup_\omega \xi$ is an upper bound of the $L^2$ convergence factor \cite{clement_mini_02_Thery_2021} and $\xi_0$ is an approximation of $\sup_\omega \xi$.
Both $\xi_0$ and $\xi_0^q$ are close to the convergence rate, with the exception of
$\xi_0^q$ that predicts much faster convergence than observed when $\theta=1.5$. This shows that the maximum of the convergence factor is not reached when $(\omega+f) \rightarrow 0$ in this case. Figure \ref{clement_mini_02_fig:evolutionErrNonlinear} confirms the results of Sec. \ref{clement_mini_02_sec:conv-lin}:
when considering $\alpha=\alpha_c$ constant, the fastest convergence is achieved when $\theta$ is close to 1, similarly to the DNWR algorithm.
However this does not translate into the nonlinear case, which converges faster with $\theta=1.5$.
Figure \ref{clement_mini_02_fig:robustesseEvolutionErrNonlinear}
shows that the convergence behavior with the linearized transmission condition is
similar to the nonlinear case.
As expected the convergence is faster for $\theta=1.5$
than for $\theta=1$.
We observed that those results are robust to changes in the values of the parameters in the range of interest.
Linearized transmission conditions are hence relevant to study
theoretically the
convergence properties of our nonlinear problem.
\section{Conclusion}
In this paper, we studied a SWR 
algorithm applied
to a simplified ocean-atmosphere problem. 
This problem considers nonlinear transmission conditions arising 
from wall laws representative of the ones used in Earth-System Models
and analogous to a quadratic friction law. We 
motivated the fact that the convergence analysis of such 
problems can only be done at a semi-discrete level in space due to the particular practical implementation of continuous interface conditions in actual climate models.  
Then we analytically studied the convergence properties in a
case with linear friction and in a case with linearized 
quadratic friction. We formulated the problem with a 
relaxation parameter $\theta$ in the transmission conditions and 
systematically assessed its impact on the convergence speed. 
%
%
%
%
For the two cases of interest, the convergence factors are derived and 
the asymptotic limits for small values of the frequency $\omega+f$ 
are given. This asymptotic limit allowed us to choose appropriate values
for the parameter $\theta$ to guarantee fast convergence of the algorithm.
%
The behavior of the algorithm for linear friction and linearized 
quadratic friction turns out to be different which leads to different
"optimal" values of $\theta$. 
%
%
Numerical experiments in the nonlinear case showed that the observed 
convergence behaves as predicted by the linearized quadratic friction
case whose thorough theoretical analysis is left for future work.

\begin{acknowledgement}
This work was supported by the French national research
agency through the ANR project %
COCOA
(grant 
ANR-16-CE01-0007). Part of this study was carried out within 
the project PROTEVS under the auspices of French 
Ministry of Defense/DGA, and led by Shom.
\end{acknowledgement}

\bibliographystyle{spmpsci}
\bibliography{clement_mini_02}

\begin{thebibliography}{1}
\providecommand{\url}[1]{{#1}}
\providecommand{\urlprefix}{URL }
\expandafter\ifx\csname urlstyle\endcsname\relax
  \providecommand{\doi}[1]{DOI~\discretionary{}{}{}#1}\else
  \providecommand{\doi}{DOI~\discretionary{}{}{}\begingroup
  \urlstyle{rm}\Url}\fi

\bibitem{clement_mini_02_kwok2013}
Gander, M., Kwok, F., Mandal, B.: {Dirichlet-Neumann and Neumann-Neumann
  waveform relaxation algorithms for parabolic problems}.
\newblock Electron. Trans. Numer. Anal. \textbf{45}, 424--456 (2016)

\bibitem{clement_mini_02_Marti_etal_2020}
Marti, O., Nguyen, S., Braconnot, P., Valcke, S., Lemari\'e, F., Blayo, E.: {A
  Schwarz iterative method to evaluate ocean--atmosphere coupling schemes:
  implementation and diagnostics in IPSL-CM6-SW-VLR}.
\newblock Geosci. Model Dev. \textbf{14}, 2959--2975 (2021)

\bibitem{clement_mini_02_mongeMultirate2021}
Meisrimel, P., Monge, A., Birken, P.: A time adaptive multirate
  {D}irichlet-{N}eumann waveform relaxation method for heterogeneous coupled
  heat equations.
\newblock preprint arXiv:2007.00410  (2020)

\bibitem{clement_mini_02_Mohammadi_etal_1998}
Mohammadi, B., Pironneau, O., Valentin, F.: Rough boundaries and wall laws.
\newblock Int. J. Numer. Methods Fluids \textbf{27}(1‐4), 169--177 (1998)

\bibitem{clement_mini_02_Pelletier_etal_2021}
Pelletier, C., Lemarié, F., Blayo, E., Bouin, M.N., Redelsperger, J.L.:
  Two-sided turbulent surface-layer parameterizations for computing air–sea
  fluxes.
\newblock {Quart. J. Roy. Meteorol. Soc.} \textbf{47}(736), 1726--1751 (2021)

\bibitem{clement_mini_02_Thery_2021}
Thery, S.: {{\'E}tude num{\'e}rique des algorithmes de couplage
  oc{\'e}an-atmosph{\`e}re avec prise en compte des param{\'e}trisations
  physiques de couches limites}.
\newblock Phd thesis, {Universit{\'e} Grenoble Alpes} (2021).
\newblock Https://tel.archives-ouvertes.fr/tel-03164786

\bibitem{clement_mini_02_Wu2017}
Wu, S.L., Al-Khaleel, M.: Optimized waveform relaxation methods for {RC}
  circuits: Discrete case.
\newblock Esaim Math. Model. Numer. Anal. \textbf{51}, 209--222 (2017)

\end{thebibliography}
\end{document}